\documentclass[final,3p,times]{elsarticle}

\usepackage{graphicx}
\usepackage{amssymb}
\usepackage{amsfonts}
\usepackage{amsmath}
\usepackage{color}
\usepackage{subfigure}
\usepackage{multirow}
\usepackage{pdflscape}

\newtheorem{theorem}{Theorem}

\newdefinition{remark}{Remark}
\newproof{pf}{Proof}
\newproof{pot}{Proof of Theorem \ref{thm2}}
\newtheorem{example}{Example}
\newdefinition{defin}{Definition}
\newdefinition{coro}{Corollary}

\begin{document}

\begin{frontmatter}

\title{A mean square chain rule and its applications in solving the random Chebyschev}

\author[l2]{J.-C. Cort\'{e}s}
\ead{jccortes@imm.upv.es}
\author[l3]{L. Villafuerte\corref{cor1}}
\ead{laura.villafuerte@unach.mx}
\author[l2]{C. Burgos}
\ead{clara.burgos@uv.es}

\cortext[cor1]{Corresponding author. Tel.: +34 (96)3879144.}

\address[l2]{Instituto Universitario de Matem\'{a}tica Multidisciplinar, Building 8G 2$^{nd}$ floor Access C, Universitat Polit\`{e}cnica de Val\`{e}ncia, 46022 Valencia (Spain)}
\address[l3]{Facultad de Ciencias en F\'{i}sica y Matem\'{a}ticas, 
Universidad Aut\'{o}noma de Chiapas, Ciudad Universitaria, 29050, 
Tuxtla Guti\'{e}rrez, Chiapas, M\'{e}xico.
Department of Mathematics,
University of Wisconsin-Madison,
Madison, Wisconsin 53706,
USA \vspace{-34pt}}

\begin{abstract}
In this paper a new version of the chain rule for calculating the mean square  derivative of a second-order stochastic process is proven.  This random operational calculus rule is applied to construct a rigorous mean square solution of the random Chebyshev differential equation (r.C.d.e.) assuming mild moment hypotheses on the random variables that appear as coefficients and initial conditions of the corresponding initial value problem. Such solution is represented through a mean square random power series. Moreover,  reliable approximations for the mean and standard deviation functions to the solution stochastic process  of the r.C.d.e. are given. Several examples, that illustrate the theoretical results, are included.
\end{abstract}
\end{frontmatter}

\section{Introduction and Preliminaries}\label{seccion 0}
Several studies of random differential equations with random coefficients have been undertaken lately, \cite{Gema10,Magdy_Sohaly_OJDM_2011, Khodabin_Maleknejad_Rostami_Nouri_MCM_2011,Santos-Dorini-Cunha-AMC-2010,Chen_Gilberto_Abraham-AMC-2014,Beltagy_Magdy_APM_2013,Knouri_MJOM_2015}. A rigorous treatment of the mean square solutions of random differential equations requires some operational tools. In \cite{Braumann_CAMWA_2010} a chain rule for the composition of a $\mathcal{C}^{1}$-function with a stochastic process  was provided. In the present paper a new version of the chain rule is proven.  While this chain rule can be applied in different scenarios,  a single application is presented in this paper, namely, the rigorous solution, in the mean square sense, of the random Chebyshev differential equation (r.C.d.e.).  

For the sake of clarity,   in the following we  summarize the main definitions and results that will be used throughout this paper. Further details about them can be found in \cite{Oksendall-LIBRO, Soong, Braumann_CAMWA_2010,Wong_Hajek-LIBRO}. 
Let $p\geq 1$ be a real number. A real random variable (r.v.) $X: \Omega \longrightarrow \mathbb{R}$ defined on a complete probability space $(\Omega, \mathcal{F},\mathbb{P})$ is called of order $p$ (in short, $p$-r.v.), if
\[
\mathbb{E}\!\left[\left|X\right|^{p}\right]<+\infty\,,
\]
where $\mathbb{E}\!\left[\,\,\right]$ denotes the expectation operator. The set $\mathrm{L}_p(\Omega)$ of all the $p$-r.v.'s endowed with the norm
\[
\left\|X \right\|_p=\left(\mathbb{E}\!\left[\left| X \right|^{p}\right]\right)^{1/p}\,,
\]
is a Banach space, \cite[p.9]{Arnold-LIBRO}. In the particular case that $X\in\mathrm{L}_2(\Omega)$, it is termed a second-order random variable (2-r.v.). Let $\mathcal{T}\subset \mathbb{R}$ be an interval of the real line, if $\mathbb{E}\!\left[\left|X(t)\right|^{p}\right]<+\infty$ for all $t\in \mathcal{T}$, then $\left\{X(t):\, t\in \mathcal{T}\right\}$ is termed a stochastic process of order $p$ (in short, $p$-s.p.).  In the particular case that $X(t)\in\mathrm{L}_2(\Omega)$ for every $t$, it is termed a second-order stochastic process (2-s.p.).

Let $\left\{X_n:\, n\geq 0  \right\}$ be a sequence in $\mathrm{L}_p(\Omega)$. We say that it is convergent in the $p$-th mean to $X\in \mathrm{L}_p(\Omega)$, if
\[
\lim_{n\rightarrow \infty} \left\| X_n-X \right\|_p=0.
\]
The so-called mean square (m.s.) convergence corresponds to $p=2$.
The  $p$-continuity and $p$-differentiability in $(\mathrm{L}_{p}(\Omega, \left\|\cdot \right\|_p)$ Banach spaces are inferred, as usually, from the $p$-norm. For instance, for the particular case  $p=2$, the s.p. $\left\{X(t):\, t\in \mathcal{T}\right\}$ in $\mathrm{L}_2(\Omega)$ is said to be mean square (m.s.)  continuous at $t\in \mathcal{T}$ if
\begin{equation}\label{p-th continuidad}
\left\| X(t+h)-X(t)\right\|_2 
\xrightarrow[h  \to 0]{} 0\,,\quad t,t+h\in \mathcal{T}.
\end{equation}
If there exists a s.p. $\frac{\mathrm{d}X(t)}{\mathrm{d}t}\in \mathrm{L}_2(\Omega)$ such that
\begin{equation}\label{p-th derivabilidad}
\left\| \frac{X(t+h)-X(t)}{h}-\frac{\mathrm{d}X(t)}{\mathrm{d}t}\right\|_2 \xrightarrow[h  \to 0]{} 0\,,\quad t,t+h\in \mathcal{T},
\end{equation}
then we say that the s.p. $X(t)$ is m.s. differentiable at $t\in \mathcal{T}$ and its  m.s. derivative at $t$ is given by $\frac{\mathrm{d}X(t)}{\mathrm{d}t}$. This calculus defined in 
$\mathrm{L}_{2}(\Omega)$ is called mean square calculus, while mean fourth (m.f.) calculus is the one associated to $p=4$ \cite{Braumann_CAMWA_2010}. In the general case that $p\geq 1$, it is usually referred to as $\mathrm{L}_{p}(\Omega)$
stochastic calculus.

In this paper we firstly establish a chain rule for calculating, in the m.s. sense, the derivative of a s.p., say $X(t)$, that is defined through the composition of another s.p., $Y(s)$, and a deterministic function, $g(t)$, i.e., $X(t)=Y(g(t))$. Afterwards, we apply such stochastic chain rule for solving the random Chebyshev differential equation. Both contributions constitute extensions of previous works by the authors and other coauthors. 

On the one hand in \cite{Braumann_CAMWA_2010} a random chain rule for differentiating, in the m.s. sense, a s.p. $X(t)=g(Y(t))$ resulting from the composition of a
deterministic function, $g(s)$, and a s.p., $Y(t)$, is established. As shown in  \cite{Braumann_CAMWA_2010}, the rigorous proof of this counterpart of the mean square stochastic chain rule needs to assume hypotheses involving mean fourth information on s.p. $Y(t)$. While   only assuming m.s. conditions of the s.p. $Y(s)$,  a random chain rule will be established    later. Therefore, a part of this paper complements the previous contribution  \cite{Braumann_CAMWA_2010}.

On the other hand, in previous contributions  important second-order differential equations have been randomized and studied taking advantage of the $\mathrm{L}_{p}(\Omega)$ stochastic calculus with $p=2,4$. In this regard, random Airy, Hermite, Legendre, Laguerre and Bessel differential equations have been studied in \cite{Airy2010,Golmankhaneh_RRP_2013,Calbo_AMC_2011,Legendre2011,Laguerre2015,Bessel_2016}, respectively. In this paper, we enlarge this list studying the randomized Chebyshev differential equation using a completely different approach based on a new m.s. chain rule.

This paper is organized as follows. In Section \ref{regla_cadena} a version of the classical chain rule is extended to compute the derivative, in the mean square sense,  of a stochastic process resulting from the composition of a stochastic process and a deterministic function. In Section \ref{rCde}, the new  stochastic chain rule is applied to solve, rigorously, the random Chebyshev differential equation by assuming quite mild assumptions on input data. Section \ref{media_varianza} is addressed to compute the main statistical properties of the solution stochastic process to the random Chebyshev differential equation, namely, the mean and the standard deviation functions, as well as providing some illustrative examples where these statistical functions are computed. The numerical results are compared with the ones obtained by Monte Carlo sampling. Our approach demonstrates its superiority against Monte Carlo simulations regarding computational timing. Finally, conclusions are drawn in Section \ref{conclusiones}.

\section{A mean square chain rule}\label{regla_cadena}

Now we  prove  a  m.s. chain rule for  the composition of a deterministic differentiable function with a s.p.  This result complements the one established in  \cite{Braumann_CAMWA_2010}.
\begin{theorem}\label{T1}
Let  $g$ be a deterministic continuous function on  $[a,b]$ such that $\frac{\mathrm{d}g(t)}{\mathrm{d}t}$ exists and is finite at some point $t\in [a,b]$. If $\{Y(s): s\in I\}$ is  a 2-s.p. such that
\begin{itemize}
\item[i)] The interval $I$ contains the range of $g$, $g([a,b])\subset I$.
\item[ii)] $Y(s)$ is m.s. differentiable at the point $g(t)$.
\item[iii)]The m.s. derivative of $Y(s)$, $\frac{\mathrm{d} Y}{\mathrm{d} s}$,  is m.s. continuous on $I$.
\end{itemize}
Then, the 2-s.p. $Y(g(t))$ is m.s. differentiable at $t$ and the m.s. derivative is given by
\[
\frac{\mathrm{d}Y(g(t))}{\mathrm{d}t}=\frac{\mathrm{d}Y}{\mathrm{d}s}\Big\rvert_{s=g(t)}\frac{\mathrm{d}g(t)}{\mathrm{d}t}.
\]
\end{theorem}
\textbf{Proof.} Let us denote $\frac{\mathrm{d}Y}{\mathrm{d}s}\Big\rvert_{s=g(t)}=\frac{\mathrm{d} Y}{\mathrm{d} s}\left(g(t)\right)$. Setting
\begin{equation}\label{r0}
G(t,\Delta t):=\frac{g(t+\Delta t)-g(t)}{\Delta t} -\frac{\mathrm{d}g(t)}{\mathrm{d}t},
\end{equation}
it follows, by hypothesis, that 
\[
 G(t,\Delta t) \xrightarrow[\Delta t \to 0]{} 0.
\]
Isolating $g(t+\Delta t)$ from (\ref{r0}), one obtains
\[
g(t+\Delta t)=g(t)+\Delta t\left(G(t,\Delta t)+\frac{\mathrm{d}g(t)}{\mathrm{d}t}\right).
\]
Let us write 
\begin{equation}\label{def_K}
K=K(t,\Delta t)=g(t+\Delta t)-g(t),
\end{equation}
and 
\begin{equation}\label{r3}
S(t,\Delta t)=
\left\{
	\begin{array}{ccl}
		\dfrac{Y(g(t+\Delta t))-Y(g(t))}{K} - \dfrac{\mathrm{d} Y}{\mathrm{d} s}\left(g(t)\right) & \mbox{if} & K \neq 0, \\[0.3cm]
		0 & \mbox{if} & K= 0.
	\end{array}
\right.
\end{equation}
Then, defining $T(t,\Delta t)=\frac{K(t,\Delta t)}{\Delta t}$, $\Delta t\neq 0$ and observing from \eqref{def_K} that $Y(g(t+\Delta t))=Y(g(t)+K)$,  expression (\ref{r3}) leads to
\begin{equation}\label{r4}
\frac{Y(g(t+\Delta t))-Y(g(t))}{\Delta t}=T(t, \Delta t)\left(S(t, \Delta t)+\frac{\mathrm{d} Y}{\mathrm{d} s}\left(g(t)\right)\right) .
\end{equation}
In view of (\ref{r4}), we must prove
\begin{equation}\label{r5}
\lim_{\Delta t \to 0}\left\| T(t, \Delta t)\left(S(t, \Delta t)+\frac{\mathrm{d} Y}{\mathrm{d} s}\left(g(t)\right) \right)-\frac{\mathrm{d} Y}{\mathrm{d} s}\left(g(t)\right) \frac{\mathrm{d}g(t)}{\mathrm{d}t} \right\|_{2}=0.
\end{equation}
From (\ref{r0}), it follows that $T(t, \Delta t)$ can be written as $T(t, \Delta t)=G(t, \Delta t)+\frac{\mathrm{d}g(t)}{\mathrm{d}t}$. Thus, by the triangle inequality one gets

\begin{eqnarray}\label{DI}
\begin{array}{c}
\left\| T(t, \Delta t)\left(S(t, \Delta t)+\frac{\mathrm{d} Y}{\mathrm{d} s}\left(g(t)\right) \right)-\frac{\mathrm{d} Y}{\mathrm{d} s}\left(g(t)\right) \frac{\mathrm{d}g(t)}{\mathrm{d}t} \right\|_{2} \\ [0.3cm] 
=\left\| G(t, \Delta t)S(t, \Delta t)+\frac{\mathrm{d}g(t)}{\mathrm{d}t}S(t, \Delta t)+G(t, \Delta t)\frac{\mathrm{d} Y}{\mathrm{d} s}\left(g(t)\right)\right\|_{2} \\  [0.3cm] 
\leq \left \| G(t, \Delta t)S(t, \Delta t)\right\|_{2}+\left\|\frac{\mathrm{d}g(t)}{\mathrm{d}t}S(t, \Delta t)\right\|_{2}+\left\|G(t, \Delta t)\frac{\mathrm{d} Y}{\mathrm{d} s}\left(g(t)\right)\right\|_{2} \\  [0.3cm] 
= \left | G(t, \Delta t)\right |
\left\|S(t, \Delta t)\right\|_{2}+\left |\frac{\mathrm{d}g(t)}{\mathrm{d}t}\right|
\left\|S(t, \Delta t)\right\|_{2}+\left |G(t, \Delta t)\right |
\left\| \frac{\mathrm{d} Y}{\mathrm{d} s}\left(g(t)\right)\right\|_{2},
\end{array}
\end{eqnarray}
where in the last step we have used that $g(t)$, and hence $\frac{\mathrm{d}g(t)}{\mathrm{d}t}$ and $G(t,\Delta t)$, are deterministic.

As $\frac{\mathrm{d} Y}{\mathrm{d} s}$ is m.s. continuous on $I$,  it is m.s. Riemann integrable too, hence by the fundamental theorem of the m.s. calculus (see \cite[p. 104]{Soong}), it follows that
\begin{equation}
\label{FCT}
Y(g(t+\Delta t)\vee g(t))-Y(g(t+\Delta t)\wedge g(t))=\int_{g(t+\Delta t)\wedge g(t)}^{g(t+\Delta t)\vee g(t)}\frac{\mathrm{d} Y(s)}{\mathrm{d} s}\,\mathrm{d}s,
\end{equation}
where $\vee$ and $\wedge$ stand for maximum and minimum operators, respectively. Using (\ref{FCT}) and the definition of  $S(t, \Delta t)$, given in (\ref{r3}), one gets
\begin{equation}\label{r6}
\left\|S(t, \Delta t)\right\|_2=
\left\{
	\begin{array}{ccl}
		 \dfrac{\displaystyle \left\| \int_{g(t+\Delta t)\wedge g(t)}^{g(t+\Delta t)\vee g(t)} \left[\frac{\mathrm{d} Y(s)}{\mathrm{d} s}-\frac{\mathrm{d} Y}{\mathrm{d} s}\left(g(t)\right) \right]\, \mathrm{d}s \right\|_2}{|K|}  & \mbox{if} & K \neq 0, \\
		0 & \mbox{if} & K= 0,
	\end{array}
\right.
\end{equation}
Moreover, as $\frac{\mathrm{d} Y}{\mathrm{d} s}$ is m.s. continuous on $I$, then by property (3) of \cite[p. 102]{Soong} we have
\begin{equation}\label{r7}
\begin{split}
&\left\| \int_{g(t+\Delta t)\wedge g(t)}^{g(t+\Delta t)\vee g(t)}\left[ \frac{\mathrm{d} Y(s)}{\mathrm{d} s}-\frac{\mathrm{d} Y}{\mathrm{d} s}\left(g(t)\right) \right] \mathrm{d}s\right\|_{2}\\
&\leq \left| K \right|\max_{s\in[g(t+\Delta t)\wedge g(t),g(t+\Delta t)\vee g(t)]}\left\|\frac{\mathrm{d} Y(s)}{\mathrm{d} s}-\frac{\mathrm{d} Y}{\mathrm{d} s}\left(g(t)\right)\right\|_{2}.
\end{split}
\end{equation}
Taking into account the continuity of $g$ on $[a,b]$ and the m.s. continuity of $\frac{\mathrm{d} Y}{\mathrm{d} s}$ on $I$, from (\ref{r6})--(\ref{r7}) one concludes that
\begin{equation}
\label{Slimit}
\lim_{\Delta t\to 0} \left \| S(t, \Delta t)\right \|_{2}=0.
\end{equation}
As $\left\| \frac{\mathrm{d} Y}{\mathrm{d} s}\left(g(t)\right) \right\|_{2}<+\infty$, $\left |\frac{\mathrm{d}g(t)}{\mathrm{d}t}\right|<+\infty$,  $G(t,\Delta t) \xrightarrow[\Delta t \to 0]{} 0$
 and $\lim_{\Delta t\to 0} \left \| S(t, \Delta t)\right \|_{2}=0$, from (\ref{DI}) it follows (\ref{r5}). Thus, the proof is complete. $\square$
\begin{example}\label{ej1}
Let $Y(s)$ be  a 2-s.p. twice m.s. differentiable such that its  second m.s. derivative, $\frac{\mathrm{d}^{2} Y}{\mathrm{d}s^2}$, exists and is m.s. continuous on the interval $I=(-1,1)$. Let  $g(t)=\ cos(t)$, $t\in (0,\pi)$. As $\frac{\mathrm{d} Y}{\mathrm{d} s}$ is m.s. differentiable on $I$, by Property (1) of \cite[p. 95]{Soong}, $\frac{\mathrm{d} Y}{\mathrm{d} s}$ is m.s. continuous on $I$. Moreover, for all $t \in (0,\pi)$,   $g(t)\in I$ and $Y(s)$  is m.s. differentiable at g(t).
Therefore, Theorem \ref{T1} asserts that, for every $t\in (0,\pi)$ the m.s. derivative of $Y(g(t))$ exists and is given by
\begin{equation}
\label{pder}
\frac{\mathrm{d}Y(g(t))}{\mathrm{d}t}=\frac{\mathrm{d} Y}{\mathrm{d} s}\Big\rvert_{s=g(t)}\frac{\mathrm{d}g(t)}{\mathrm{d}t}=\frac{\mathrm{d} Y}{\mathrm{d} s}\left(\cos(t)\right)\left[-\sin(t)\right].
\end{equation}
Now, taking into account that $\frac{\mathrm{d}^{2} Y}{\mathrm{d}s^2}$ is m.s. continuous on $I$, with the aid of Property (4) of \cite[p.96]{Soong} and Theorem \ref{T1}, it follows that for every $t\in (0,\pi)$, the second m.s. derivative of $Y(g(t))$ is
\begin{equation*}
\frac{\mathrm{d}^{2}Y(g(t))}{\mathrm{d}t^{2}}=\frac{\mathrm{d}\left(\frac{\mathrm{d} Y}{\mathrm{d} s}\left(\cos(t)\right)\left[-\sin(t)\right]\right)}{\mathrm{d}t}=\sin^{2}(t)\frac{\mathrm{d}^{2} Y}{\mathrm{d}s^2}(\cos (t) )-\cos (t)\frac{\mathrm{d}Y}{\mathrm{d}s}(\cos (t)).
\end{equation*}
Using (\ref{pder}), $\frac{d^{2}Y(g(t))}{dt^{2}}$ can be written as 
\begin{equation}
\label{sder}
\frac{\mathrm{d}^{2}Y(g(t))}{\mathrm{d}t^{2}}=\sin^{2}(t)\frac{\mathrm{d}^{2} Y}{\mathrm{d}s^2}(\cos (t))+\frac{\cos(t)}{\sin(t)}\frac{\mathrm{d}Y(g(t))}{\mathrm{d}t}.
\end{equation}
\end{example}
\section{Solving the random Chebyshev differential equation (r.C.d.e.)}\label{rCde}
Let $A, Y_{0}$ and $Y_{1}$ be r.v.'s and $I=(-1,1)$.  Consider  the r.C.d.e. with two initial conditions:
\begin{align}
\label{RCDE}
&(1-s^{2})\frac{\mathrm{d}^2Y(s)}{\mathrm{d}s^2}-s\frac{\mathrm{d}Y(s)}{\mathrm{d}s}+A^{2}Y(s)=0, \quad s\in I,\\
&Y(0)=Y_{0},\qquad \frac{\mathrm{d}Y}{\mathrm{d}s}(0)=Y_{1}.\label{IV}
\end{align}
In the following, the   chain rule that has been established in Theorem \ref{T1},  will be applied to find a 2-s.p., $\{Y(s): s\in I\}$, satisfying  the initial value problem (\ref{RCDE})--(\ref{IV}). As usually, the 2-s.p. $Y(s)$ is called the m.s. solution of (\ref{RCDE})--(\ref{IV}). Let us assume that the second m.s. derivative of $Y(s)$, $\frac{\mathrm{d}^2Y}{\mathrm{d}s^2}$,  exists and  is m.s.  continuous on $I$.
Setting $s=\cos(t):=g(t)$, $t\in (0,\pi)$, Theorem \ref{T1} (chain rule) and 
Example \ref{ej1} (see expressions (\ref{pder})--(\ref{sder})) yield
\begin{equation}\label{r10}
\frac{\mathrm{d}Y}{\mathrm{d}s}(s)=\frac{\mathrm{d}Y}{\mathrm{d}s}(\cos(t))=-\frac{1}{\sin(t)}\frac{\mathrm{d}Y(\cos(t))}{\mathrm{d}t},
\end{equation}
and 
\begin{equation}\label{r11}
\frac{\mathrm{d}^2Y}{\mathrm{d}s^2}(s)=\frac{\mathrm{d}^2Y}{\mathrm{d}s^2}(\cos(t))=\frac{1}{\sin^{2}(t)}\left(\frac{\mathrm{d}^2Y(\cos(t))}{\mathrm{d}t^2}
-\frac{\cos(t)}{\sin(t)}\frac{\mathrm{d}Y(\cos(t))}{\mathrm{d}t}\right).
\end{equation}
Defining 
\begin{equation}\label{XvsY}
X(t)=Y(\cos(t)),\quad s=\cos(t),
\end{equation}
and inserting (\ref{r10}) and (\ref{r11}) in expression (\ref{RCDE}) we have
\[
(1-\cos^{2}(t))\frac{1}{\sin^{2}(t)}\left[\frac{\mathrm{d}^{2}X(t)}{\mathrm{d}t^{2}}-\frac{\cos(t)}{\sin(t)}\frac{\mathrm{d}X(t)}{\mathrm{d}t}\right]-\cos(t)\left[-\frac{1}{\sin(t)}\frac{\mathrm{d}X(t)}{\mathrm{d}t}\right]+A^{2}X(t)=0.
\]
Hence, using the well-known identity $\cos^{2}(t)+\sin^{2}(t)=1$, one gets
\[
\frac{\mathrm{d}^{2}X(t)}{\mathrm{d}t^{2}}+A^{2}X(t)=0,\qquad t\in(0,\pi).
\]
Now, if $\cos(t)=0$, $t\in(0,\pi)$ then $t=\frac{\pi}{2}$. Therefore, by \eqref{XvsY}, one gets
\[
Y_{0}:=Y(0)=Y(\cos(\pi/2))=X(\pi/2).
\]
Also, by \eqref{r10} and \eqref{XvsY} one obtains

\begin{equation*}
Y_{1}:=\frac{\mathrm{d}Y}{\mathrm{d}s}(0)=\frac{\mathrm{d}Y}{\mathrm{d}s}(\cos(\pi/2))=-\frac{1}{\sin(\pi/2)} \frac{\mathrm{d}X}{\mathrm{d}t}\Big\rvert_{t=\pi/2}=-\frac{\mathrm{d}X}{\mathrm{d}t}(\pi/2).
\end{equation*}

Therefore, the change of variable $s=\cos(t)$ transforms the problem (\ref{RCDE})--(\ref{IV}) into the random initial value problem
\begin{align} \label{1RCDE}
    &\frac{\mathrm{d}^{2}X(t)}{\mathrm{d}t^{2}}+A^{2}X(t)=0,\qquad t\in(0,\pi),   \\ 
    & X(\pi/2)=Y_{0},\qquad  \frac{\mathrm{d}X}{\mathrm{d}t}(\pi/2)=-Y_{1}. \label{XRCDE}
\end{align}
Moreover, letting $r=t-\pi/2$ and defining $Z(r):=X(r+\pi/2)=X(t)$ it follows $$\frac{\mathrm{d}Z(r)}{\mathrm{d}r}=\frac{\mathrm{d}X(t)}{\mathrm{d}t},\qquad \frac{\mathrm{d}^{2}Z(r)}{\mathrm{d}r^{2}}=\frac{\mathrm{d}^{2}X(t)}{\mathrm{d}t^{2}},$$
and, according to \eqref{XRCDE} one gets
$$Z(0)=X(\pi/2)=Y_{0},\qquad  \frac{\mathrm{d}Z}{\mathrm{d}r}(0)= \frac{\mathrm{d}X}{\mathrm{d}t}(\pi/2)=-Y_{1}.$$ Finally, the problem  (\ref{RCDE})--(\ref{IV}) has been transformed into
\begin{align}
\label{2RCDE}
    &\frac{\mathrm{d}^{2}Z(r)}{\mathrm{d}r^{2}}+A^{2}Z(r)=0,\qquad r\in(-\pi/2,\pi/2),   \\
    & Z(0)=Y_{0},\qquad  \frac{\mathrm{d}Z}{\mathrm{d}r}(0)=-Y_{1}.\label{2IV}
\end{align}
In the following result a m.s. solution of the r.C.d.e. is provided. Its proof is a direct application of Theorem 3.1 of \cite{Gema10}.

\begin{theorem}\label{T2}
Let $A^{2}$  be a r.v. such that there exist  constants $M>0$ and $\kappa$ such that  $0\leq \kappa <2$, satisfying the following property
\begin{equation}\label{RP}
\left\| (A^{2})^{n}\right\|_{4}=\mathcal{O}(M^{n-1}((n-1)!)^{\kappa}), \qquad \forall n\geq 1.
\end{equation} 
If $Y_{0}$ and $Y_{1}$  are 4-r.v.'s both independent of r.v. $A^{2}$. Then, a m.s. solution of (\ref{RCDE})--(\ref{IV})  is given by:

\begin{equation}\label{3sol}
\begin{split}
Y(s)&=\sum_{k\geq 0}\frac{(-1)^{k}(A^{2})^{k}Y_{0}}{(2k)!}(\cos^{-1}(s)-\pi/2)^{2k}\\ 
&+\sum_{k\geq 0}\frac{(-1)^{k}(A^{2})^{k}(-Y_{1})}{(2k+1)!}(\cos^{-1}(s)-\pi/2)^{2k+1},\quad
s\in(-1,1).
\end{split}
\end{equation}
\end{theorem}
\textbf{Proof.} From Theorem 3.1 of \cite{Gema10} it follows that
\[
Z(r)=\sum_{k\geq 0}\frac{(-1)^{k}(A^{2})^{k}Y_{0}}{(2k)!}r^{2k}+\sum_{k\geq 0}\frac{(-1)^{k}(A^{2})^{k}(-Y_{1})}{(2k+1)!}r^{2k+1},\quad
r\in(-\pi/2,\pi/2)
\]
is a m.s. solution of problem  (\ref{2RCDE})--(\ref{2IV}). As $X(t)=X(r+\pi/2)=Z(r)$ then 
\begin{equation}
\label{2sol}
X(t)=\sum_{k\geq 0}\frac{(-1)^{k}(A^{2})^{k}Y_{0}}{(2k)!}(t-\pi/2)^{2k}+\sum_{k\geq 0}\frac{(-1)^{k}(A^{2})^{k}(-Y_{1})}{(2k+1)!}(t-\pi/2)^{2k+1}
\end{equation}
is a m.s. solution of (\ref{1RCDE})--(\ref{XRCDE}) for $t\in (0,\pi)$. Finally, the result follows by considering relation \eqref{XvsY}.
$\square$
\begin{remark}\label{r1}
If $A$ is a r.v. such that  there exists $0\leq \kappa <2$ satisfying the property
\begin{equation}\label{RP1}
\left\| (A^{2})^{n+1}\right\|_{4}=\mathcal{O}(n^{\kappa})\left \|(A^{2})^{n}\right\|_{4},\qquad \forall n\geq 0.
\end{equation} 
Then $A^{2}$ also satisfies property (\ref{RP}), see Lemma 2.2 of \cite{Gema10}. Sometimes, verifying property (\ref{RP1}) might be easier than checking property (\ref{RP}).
\end{remark}
\begin{remark}\label{nv} In Example 2.3. of \cite{Gema10} it is shown that any beta r.v. satisfies property (\ref{RP}). In the following, we add some further important families of r.v.'s that also satisfy such property, thus enlarging and enriching the applicability of our theoretical results. Firstly, an important class of r.v.'s that satisfies property (\ref{RP}) is the centered gaussian family. Indeed, a  gaussian r.v.  $A$ with zero-mean and arbitrary standard deviation, $\sigma>0$, satisfies condition (\ref{RP1}) with $\kappa=1$:
\begin{equation}\label{momentos_normal_media_0}
\frac{\left\| (A^{2})^{n+1}\right\|_{4}}{\left\| (A^{2})^{n}\right\|_{4}}=\frac{\sigma^{2}}{2}\left(\frac{\prod_{i=1}^{8}(8n+i)}{\prod_{i=1}^{4}(4n+i)}\right)^{1/4}=\mathcal{O}(n),\qquad A\sim \mathrm{N}(0;\sigma ^{2}).
\end{equation} 
Thus, by Remark \ref{r1}, $A^{2}$ satisfies property (\ref{RP}). Notice that in \eqref{momentos_normal_media_0}, we have used the following property
$$
\mathbb{E}\left[ X^{2m} \right]
=
\dfrac{(2m)!}{2^m m!} \sigma^{2m}, \quad m=0,1,2,\ldots, \qquad X\sim \mathrm{N}(0;\sigma ^{2}),
$$
for $X=A$ and $m=4n$ and $m=4(n+1)$.  Another important class of r.v.'s satisfying (\ref{RP}) are  bounded r.v.'s. In fact, let us consider that $A$ is  an absolutely continuous  r.v. with probability density function (p.d.f.) $f_A(a)$ such that $a_1 \leq A(\omega) \leq a_2$ for every $\omega \in \Omega$. Then,   $|A^{2}(\omega)|\leq H$ being $H=\max\{1, (a_1)^2,(a_2)^2 \}$, and 
\begin{equation}\label{momentos_truncada}
\left\| (A^{2})^{n}\right\|_{4}
=
\left(
\mathbb{E}
\left[
A^{8n}
\right]
\right)^{\frac{1}{4}}
=
\displaystyle
\left(
\int_{-\infty}^{\infty}
a^{8n} f_A(a) \mathrm{d}a
\right)^{\frac{1}{4}}
\leq
H^{2n} \displaystyle
\left(
\int_{-\infty}^{\infty}
 f_A(a) \mathrm{d}a
\right)=H^{2n},
\end{equation}
where in the last step we have used that the integral is the unit since $f_A(a)$ is a p.d.f. Notice that \eqref{momentos_truncada} can equivalently be  written in the following form
$$
\left\| (A^{2})^{n}\right\|_{4} \leq L\times L^{n-1},\quad
 L=H^2.
$$
Therefore,  property  (\ref{RP}) holds for $\kappa=0$ and $M=L>0$. The previous reasoning is also valid for discrete r.v.'s simply by substituting the integral and the corresponding series.  We finally point out that this family of r.v.'s is particularly useful in applications because apart from embracing important classes of bounded standard r.v.'s, like binomial, uniform, beta, etc., any unbounded r.v. can be fairly approximated by truncating its domain  adequately. For instance, according to the probabilistic  Chebyshev's inequality, for every 2-r.v. $X$ with mean $\mu_X$ and variance $\sigma_X^2>0$, one gets
$$
\mathbb{P}\left[
|X-\mu_X|>k\sigma
\right]
\leq \dfrac{1}{k^2},\quad k>0.
$$
Hence taking $k=10$, the $99\%$ of the probability mass of any 2-r.v. is contained in the truncated domain $[\mu_X-10\sigma_X,\mu_X+10\sigma_X]$. Naturally, the diameter of this interval can considerably be reduced if the distribution of $X$ is known.
\end{remark}
\begin{remark}\label{rt}
If $A^{2}$ is a r.v. satisfying property (\ref{RP}), then the trigonometric processes
$$\sin(At)=\sum_{k\geq 0}\frac{(-1)^{k}(A^{2})^{k+\frac{1}{2}}}{(2k+1)!}t^{2k+1}, \quad 
\cos(At)=\sum_{k\geq 0}\frac{(-1)^{k}(A^{2})^{k}}{(2k)!}t^{2k},$$
are well defined in $\mathrm{L}_{2}(\Omega)$. Thus, the solution s.p. $Y(s)$ of the r.C.d.e. defined by \eqref{3sol}, can be written as
\begin{equation}
\label{scsol}
Y(s)=Y_{0}\cos \left(A \left[cos^{-1}(s)-\frac{\pi}{2} \right] \right)-\frac{Y_{1}}{A}\sin\left(A\left[\cos^{-1}(s)-\frac{\pi}{2}\right]\right).
\end{equation}
Moreover,  from Remark 3.3 of \cite{Gema10}, it follows
$$
\cos\left(A \left[\cos^{-1}(s)-\frac{\pi}{2}\right]\right)=\cos\left(A\cos^{-1}(s)\right)\cos\left(A\frac{\pi}{2}\right)+\sin\left(A\cos^{-1}(s)\right)\sin\left(A\frac{\pi}{2}\right),
$$
and
$$\sin\left(A\left[\cos^{-1}(s)-\frac{\pi}{2}\right]\right)=\sin\left(A\cos^{-1}(s)\right)\cos\left(A\frac{\pi}{2}\right)
-\sin\left(A\frac{\pi}{2}\right)\cos\left(A\cos^{-1}(s)\right).
$$
\end{remark}
\section{Statistical properties of the solution. Examples}\label{media_varianza}
A natural application of previous results is the extension to the random framework of classical Chebyshev polynomials. If $A$ is a bounded r.v. that  only takes positive integers values, then the random Chebyshev polynomials of the first kind, denoted by $T_{A}$, are defined by
 $$T_{A}(s)=\cos(A\cos^{-1}(s)),\qquad s\in (-1,1).$$
Also  the random Chebyshev polynomials of the second kind, denoted by $U_{A-1}$,  are defined by the expression
 $$\sin(A\cos^{-1}(s))=U_{A-1}(s)\sin(\cos^{-1}(s)),\qquad s\in (-1,1).$$
The deterministic Chebychev polinomials of the first and the second kind can be found in a number of references, for example in \cite[pp. 57, 61]{Agarwal}. By Remark \ref{rt}, if $A$ is a r.v. taking only positive integer values, then  the solution  $Y(s)$ of  the r.C.d.e. can be written as 
\begin{equation} \label{sp}
\begin{split}
Y(s)&=\cos\left(A\frac{\pi}{2}\right)\left[Y_{0}T_{A}(s)-\frac{Y_{1}}{A}U_{A-1}(s)\sin(\cos^{-1}(s))\right]\\
 &+\sin\left(A\frac{\pi}{2}\right)\left[Y_{0}U_{A-1}\sin(\cos^{-1}(s))+\frac{Y_{1}}{A}T_{A}(s)\right].   
\end{split}
\end{equation}
\begin{example}
Suppose that $A$ is a r.v. such that $\mathbb{P}[A=2]=\mathbb{P}[A=4]=\mathbb{P}[A=6]=1/3$. If $Y_{1}=0$ w.p. 1 (with probability 1) and $Y_{0}$ is a r.v. with $\mathbb{E}[Y_{0}]=1$ and 
$\mathbb{E}[(Y_{0})^{2}]=3/2$, then according to \eqref{sp} the solution of the r.C.d.e. is given by 
$$Y(s)= \cos\left(A \dfrac{\pi}{2} \right)Y_{0}T_{A}(s)+\sin\left(A \dfrac{\pi}{2}  \right) Y_{0}\sin(A\cos^{-1}(s)).$$
Assuming that $Y_{0}$ is independent of $A$ and taking into account that $\mathbb{E}\left[  Y_{0} \right]=1$ and that only takes the following even values $A\in \{ 2,4,6 \}$, one gets
$$
\mathbb{E}\left[\sin\left(A \dfrac{\pi}{2} \right) Y_{0}\sin(A\cos^{-1}(s))\right]
=
\mathbb{E}\left[\sin\left(A \dfrac{\pi}{2} \right)  \sin(A\cos^{-1}(s))\right]
=
0,
$$
and 
$$
\mathbb{E} \left[  \cos\left(A \dfrac{\pi}{2}\right)Y_{0}T_{A}(s) \right]
=\sum_{k=1}^{3}\cos\left(2k\frac{\pi}{2}\right)T_{2k}(s)\mathbb{P}(A=2k).
$$
Further, in view of $T_{2}(s)= 2s^{2}-1$, $T_{4}(s)=8s^{4}-8s^{2}+1$, $T_{6}(s)=32s^{6}-48s^{4}+18s^{2}-1$, it follows that
\begin{equation}
\begin{split}
\mathbb{E}\left[Y(s)\right]&=\frac{1}{3}\left(-T_{2}(s)+T_{4}(s)-T_{6}(s)\right)\\
&=\frac{1}{3}\left(-32s^{6}+56s^{4}-28s^{2}+3\right),
\end{split}
\end{equation}
and 
\begin{equation} \label{sp}
\begin{split}
\mathbb{E}\left[(Y(s))^{2}\right]&=\mathbb{E}\left[(Y_{0})^{2}\right]\mathbb{E}\left[\left(\cos \left(A\frac{\pi}{2} \right) \right)^{2} (T_{A}(s))^{2}\right]\\
 &=\frac{3}{2}\sum_{k=1}^{3}\left(\cos\left(2k\frac{\pi}{2}\right)\right)^{2}\left(T_{2k}(s)\right)^{2}\mathbb{P}(A=2k)\\
 &=\frac{1}{2}\left(T_{2}(s))^{2}+(T_{4}(s))^{2}+(T_{6}(s))^{2}\right)\\
 &=-512s^{12}+1536s^{10}-1696s^{8}+832s^{6}-168s^{4}+8s^{2}+\frac{1}{2} \,.
 \end{split}
\end{equation}
\end{example}
Now, if $A$ does not take any integer values w.p. 1, but it satisfies property (\ref{RP}), then we compute approximations of the mean and standard deviation by truncating the infinite series given by (\ref{3sol}):
\begin{equation}\label{truncation}
Y_{N}(s)=Y_{0}F(s;A,N)\\ 
+Y_{1}G(s;A,N),\quad
s\in(-1,1).
\end{equation}
where 
\begin{equation}\label{series_truncadas}
\begin{array}{ccl}
F(s;A,N) &=& \displaystyle \sum_{k =0}^{N}\frac{(-1)^{k}(A^{2})^{k}}{(2k)!}(\cos^{-1}(s)-\pi/2)^{2k},\\[0.3cm]
G(s;A,N) &=& \displaystyle \sum_{k=0}^{N}\frac{(-1)^{k+1}(A^{2})^{k}}{(2k+1)!}(\cos^{-1}(s)-\pi/2)^{2k+1}.
\end{array}
\end{equation}
Under the  hypothesis that r.v. $A$ is independent of both r.v.'s $Y_{0}$ and $Y_{1}$, one obtains
\begin{equation}
\label{amean}
\mathbb{E}\left[Y_{N}(s)\right]=\mathbb{E}\left[Y_{0}\right]\mathbb{E}\left[F(s;A,N)\right]+\mathbb{E}\left[Y_{1}\right]\mathbb{E}\left[G(s;N,A)\right]
\end{equation}
and
\begin{equation}
\label{asecodm}
\begin{split}
\mathbb{E} \left[ (Y_{N})^{2} \right]&=\mathbb{E} \left[(Y_{0})^{2}\right] \mathbb{E}\left[(F(s;N,A))^{2}\right]+\mathbb{E}\left[(Y_{1}^{2})\right]\mathbb{E}\left[(G(s;N,A))^{2}\right]\\
 &+2\mathbb{E}\left[Y_{1}\right]\mathbb{E}\left[Y_{0}\right]\mathbb{E}\left[F(s;N,A)G(s;N,A)\right].
\end{split}
\end{equation}
The usefulness of  expressions \eqref{truncation}--\eqref{series_truncadas} to compute the approximations of the mean and the variance functions of the solution s.p., $Y(s)$,  according to \eqref{amean} and \eqref{asecodm}, resorts in the following key property: If $\{ Y_N:\, N\geq 0 \}$ is a sequence of r.v.'s in $\mathrm{L}_2(\Omega)$ that m.s. converges to the r.v. $Y\in\mathrm{L}_2(\Omega)$, then
\begin{equation}\label{propiedad_m_s}
\mathbb{E}[Y_N] \xrightarrow[N  \to \infty]{} \mathbb{E}[Y],\quad
 \mathbb{E}[(Y_N)^2] \xrightarrow[N  \to \infty]{} \mathbb{E}[Y^2],
\end{equation}  
and as a consequence, an analogous property for the variance
\begin{equation}\label{propiedad_m_s_2}
\sigma^2_{Y_N} \xrightarrow[N  \to \infty]{} \sigma^2_{Y}
\end{equation}
holds (see Theorem 4.3.1 in \cite[p. 88]{Soong})).

\begin{example}\label{eje2} In this second example, we consider  $A$ with normal distribution, $A\sim \mathrm{N}(0,1/4)$, $Y_{0}$ with Beta distribution,  $Y_{0}\sim \mathrm{Be}(1,3)$ and 
$Y_{1}$ with uniform distribution on $[0,2]$, $Y_{1}\sim \mathrm{U}(0,2)$. Expressions (\ref{amean}) and (\ref{asecodm}) are used to compute the mean $\mathbb{E}[Y_{N}(s)]$ and the standard deviation $\sigma_{N}(s)=\sqrt{\mathbb{E}[Y_{N}(s)^{2}]-(\mathbb{E}[Y_{N}(s)])^{2}}$. We refer to this method as the truncated series method (T.S.M.). Numerical results obtained by this approach are compared with those reached from the Monte Carlo Method (M.C.M.). In this latter case we use the notations $\tilde{\mu}_{Y}^{m}(s)$ and $\tilde{\sigma}_{Y}^{m}(t)$ for the mean and standard deviation, respectively, being $m$  the number of simulations. From  (\ref{sp})  the exact values of the mean, $\mathbb{E}[Y(s)]$, and  standard deviation, $\sigma[Y(s)]$, can be computed. In the following, this method is termed the theoretical method. According to Remark \ref{nv},  $A$ satisfies condition (\ref{RP}), therefore  $Y_{N}(s)$ is m.s. convergent in $\mathrm{L}_{2}(\Omega)$. Therefore, using  \eqref{propiedad_m_s}--\eqref{propiedad_m_s_2}, it follows that  sequences of $\mathbb{E}[Y_{N}(s)]$ and $\sigma_{Y_N}(s)$ are convergent to the corresponding exact values.  The results collected in Tables \ref{ta1}-\ref{ta2}   illustrate this assertion.  A comparison of the CPU time used in Mathematica$^{\textregistered}$ 10.3  to compute the numerical results  shown in Tables \ref{ta1}-\ref{ta2} is given in Table \ref{ta3}. Data shows that the proposed truncated series method is faster than both Monte Carlo  and the Theoretical methods.

%
%
%
%
%
%
\begin{table}[h]
   \scriptsize{ \centering 
\caption{Approximations of the mean by the proposed truncated series method ($\mathbb{E}[Y_{N}(s)]$) and Monte Carlo sampling ($\tilde{\mu}_{Y}^{m}(s)$) using $N=10$ as  truncation order,  different number of  simulations, $m$, at some selected time points $s$ in the context of Example \ref{eje2}. Last column contains the values of the exact  expactation $\mathbb{E}[Y(s)]$.}\label{ta1}
\vspace{0.1cm}
\begin{tabular}{*{5}{c}}
\hline
$s$ &  $\mathbb{E}[Y_{N}(s)]$; $N=10$  & $\tilde{\mu}_{Y}^{m}(s)$; $m=50000$& $\tilde{\mu}_{Y}^{m}(s)$; $m=100000$& $\mathbb{E}[Y(s)]$\\
\hline
 0.1&0.349812&0.349456&0.350778 &0.349812\\

0.3&0.550634& 0.550906&0.551887&0.550634\\

0.5&0.759256& 0.760167&0.760765&0.759256\\

0.7 &0.988303& 0.989880&0.990006&0.988303 \\

0.9 &1.277650 &1.279990&1.279390&1.277650\\
\hline
\end{tabular}}
\end{table}

\begin{table}[h]
   \scriptsize{ \centering 
  \caption{Approximations of the standard deviation by the proposed truncated series method ($\sigma_N(s)$) and Monte Carlo sampling ($\tilde{\sigma}_{Y}^{m}(s)$) using $N=10$ as truncation order and different  number  $m$ of simulations, at some selected time points $s$ in the context of Example \ref{eje2}. Last column contains the values of the exact  standard deviation $\sigma[Y(s)]$.}\label{ta2}
\vspace{0.1cm}
\begin{tabular}{*{5}{c}}
\hline
$s$ &  $\sigma_N(s)$; $N=10$ & $\tilde{\sigma}_{Y}^{m}(s)$; $m=50000$& $\tilde{\sigma}_{Y}^{m}(s)$; $m=100000$ & $\sigma[Y(s)]$\\
\hline

0.1&0.201862&0.202049&0.202471&0.201862\\

0.3&0.259597&0.259382&0.259939&0.259597\\

0.5&0.353343&0.352766&0.353362&0.353343\\

0.7 &0.475575&0.474643&0.475218&0.475575\\

0.9 &0.650469&0.640920&0.649552&0.649552\\

\hline
\end{tabular}}
\end{table}

\end{example}

\begin{table}[h]
   \scriptsize{ \centering 
  \caption{ Execution time  for computing the mean and variance for Example \ref{eje2} implemented on Intel$\scriptsize{^{\textregistered}}$ Core{\scriptsize$^{\texttrademark}$} 2 Duo, 4GB, 2.4GHz using the software Mathematica 10.3.}\label{ta3}
\vspace{0.1cm}
\begin{tabular}{*{4}{c}}
\hline
$Methods$ & Time  CPU (seconds) & \% Increase \\
\hline
Truncated series method (T.S.M.) && \\
truncation order $N=10$ &  0.185919&\\
\hline
Monte Carlo (M.C.) & & From  T.S.M. to M.C.  \\
$10\times 10^4$ simulations  &  103.014& 55408.3\\
\hline
Theoretical method (T.M.) && From T.S.M. to T.M. \\
&1.53099 & 822.393\\
\hline
 \end{tabular}}
\end{table}

\section{Conclusions}\label{conclusiones}
The most important result established in this paper is a chain rule for differentiating stochastic processes, in the mean square sense. This chain rule applies to  stochastic processes, say $X(t)=Y(g(t))$, where  $Y(s)$ is a mean square differentiable stochastic process and $g(t)$ is a  deterministic mapping.  This result complements  the stochastic chain rule for differentiating $X(t)=g(Y(t))$ in the mean square sense, that has  been established in \cite{Braumann_CAMWA_2010}.  It is remarkable that in this latter contribution $\mathrm{L}_p$-random calculus with $p=2$ (mean square calculus) and $p=4$ (mean fourth calculus) are required, whereas here mean square calculus is only used. As a significant application of this new stochastic chain rule, we have solved the Chebyshev random differential equation. Finally, it is important to point out that the availability of this novel stochastic chain rule opens up numerous  potential applications including the rigorous solution of other random differential equations.   

\section*{Acknowledgements}
This work has been partially supported by the Ministerio de Econom\'{i}a y Competitividad grant MTM2013-41765-P and Mexican Conacyt. This work was completed while the second author was a Visiting Fellow at UW Department of Mathematics in Madison, USA. The hospitality and support provided by this appointment is gratefully acknowledged.


\end{document}